\documentclass[a4paper, 12pt]{article}
\usepackage[utf8]{inputenc}
\usepackage[english]{babel}
\usepackage{graphics,graphicx}
\usepackage{amssymb,amsmath}
\usepackage{color}
\usepackage{enumitem}
\usepackage{authblk}
\usepackage{fullpage}
\usepackage{amsthm}
\usepackage{tikz}
\usepackage{caption}
\usepackage{subcaption}
\usepackage{todonotes}
\usepackage{hyperref, url, enumitem}

\usepackage{algorithm}%
\usepackage{algorithmicx}%
\usepackage{algpseudocode}
\usepackage{algorithm, algpseudocode}

\usetikzlibrary{arrows}
\usepackage[final]{pdfpages}

\newtheorem{theorem}{Theorem}
\newtheorem{lemma}[theorem]{Lemma}

\newtheorem{corollary}[theorem]{Corollary}

\newtheorem{observation}[theorem]{Observation}
\newtheorem{claim}{Claim}

\newtheorem{conj}[theorem]{Conjecture}
\newtheorem{definition}[theorem]{Definition}

\setlist[itemize]{leftmargin=1.4cm}
\setlist[enumerate]{leftmargin=*}
\title{Structure of (bull, diamond)-free graphs and its applications}
\author{Suchismita Mishra}
\affil{Departamento de Matem\'aticas, Univesidad Andr\'es Bello, Chile.}
\date{ }

\begin{document}
	\maketitle

 \begin{abstract}
     This paper discusses the complete structure of the (bull, diamond)-free graphs. As an application of that, we give the characterization of the partitionable (bull, diamond)-free graphs. Moreover, we show that such a partition for a partitionable (bull, diamond)-free graph can be found in polynomial time. Additionally, we show that the cop number of a (bull, diamond)-free graph containing a triangle is at most two less than its diameter. Furthermore, the cop number of a connected ($P_n$, bull, diamond)-free graph with a triangle, is at most $n-3$, for any natural number $n>3$. We also discuss a couple of applications of the structural theorem of the (bull, diamond)-free graphs in the conclusions. 
\end{abstract}


\section{Introduction}

All the graphs in this paper are simple, undirected and finite. A \emph{clique} is a set of mutually adjacent vertices. The complete graph, cycle and path on $n$ vertices are denoted by $K_n, C_n$ and $P_n$, respectively. A cycle on three vertices is called a \emph{triangle}. A \emph{hole} is a cycle of length at least four and an \emph{odd-hole} is a hole on an odd number of vertices. Let $\mathcal{F}$ be a family of graphs. A graph $G$ is said to be $\mathcal{F}$-free if no induced subgraph of $G$ is isomorphic to a graph in $\mathcal{F}$. When $\mathcal{F}$ is a finite set $\{H_1,H_2, \dots H_t\}$ then instead of writing $\mathbf{F}$-free, we write $(H_1, H_2, \dots H_t)$-free. Moreover, when the family consists of a single graph $H$ then we write $H$-free.

In the structural graph theory, most of the problems revolve around the structural behavior of the graph. In this paper, we discuss the structural behavior of graphs from a specific class. A \emph{bull} is the graph obtained from a triangle by adding two pendent edges in two vertices. In \cite{chudnovskyBullFreePartI} and \cite{chudnovskyBullFreePartIiIii}, Chudnovsky studied the structure of the bull-free graphs, which was later used to answer several questions on bull-free graphs, for example, see \cite{thomasse2017polynomial}. Another important class of graphs is the diamond-free graphs. A \emph{diamond} is the graph obtained from the complete graph on four vertices by deleting an edge. In this paper, we analyze the detailed structural behavior of all (bull, diamond)-free graphs. Then we discuss a couple of applications of the same: partitionable (bull, diamond)-free graphs, cops and robber games in (bull, diamond)-free graphs. 

A graph $G$ is called \emph{partitionable} if the vertex set can be partitioned into $(A, B)$ such that $A$ induces a disjoint union of cliques and $B$ induces a triangle-free graph. The concept of partitionable graphs is recently introduced by Abu-Khzam, Feghali and Müller \cite{KFM2015}. 
The computational complexity of deciding whether a graph is partitionable is $NP$-complete even when restricted to the class of bull-free graphs \cite{KFM2015}. Bougeret and Ochem \cite{BO2017} showed that the decision problem remains $NP$-complete for the following sub-classes of bull-free graphs.
\begin{itemize}
    \item planar ($C_4, . . . , C_t$, bull, gem, odd hole)-free graphs with maximum degree $8$.
    \item ($K_4$, bull, house, $C_5, \dots C_t$)-free graph.
\end{itemize}

\begin{figure}[h!]
  \centering
  \begin{tikzpicture}[scale=1]
   \tikzstyle{vertex}=[circle,fill,draw,thick,inner sep=0.5mm,minimum size=1.5mm]

   \node[] (u1) at (0,-1) [vertex] {};
   \node[] (u2) at (-0.5,0) [vertex] {};
   \node[] (u3) at (0.5,0) [vertex] {};
   \node[] (u4) at (-0.5,1) [vertex] {};
   \node[] (u5) at (0.5,1) [vertex] {};
   \node[] at (0,-2) {bull};

   \node[] (v1) at (2,-1) [vertex] {};
   \node[] (v2) at (1.5,0) [vertex] {};
   \node[] (v3) at (2,1) [vertex] {};
   \node[] (v4) at (2.5,0) [vertex] {};
   \node[] at (2,-2) {diamond};

   \node[] (x1) at (3.5,-0.5) [vertex] {};
   \node[] (x2) at (4,-1) [vertex] {};
   \node[] (x3) at (5,-1) [vertex] {};
   \node[] (x4) at (5.5,-0.5) [vertex] {};
   \node[] (x5) at (4.5,0.5) [vertex] {};
   \node[] at (4.5,-2) {gem};
   
   \node[] (y1) at (6.5,0) [vertex] {};
   \node[] (y2) at (6.5,-1) [vertex] {};
   \node[] (y3) at (7.5,-1) [vertex] {};
   \node[] (y4) at (7.5,0) [vertex] {};
   \node[] (y5) at (7,1) [vertex] {};
   \node[] at (7,-2) {house};

   \node[] (z1) at (8.5,1) [vertex] {};
   \node[] (z2) at (9.5,1) [vertex] {};
   \node[] (z3) at (9,0) [vertex] {};
   \node[] (z4) at (9,-1) [vertex] {};
   \node[] at (9,-2) {paw};
   

   \draw[thick] (u1)--(u2)--(u3)--(u1);
   \draw[thick] (u2) -- (u4);
   \draw[thick] (u3) -- (u5);
   \draw[thick] (v2) -- (v4);
   \draw[thick] (v1)--(v2)--(v3)--(v4)--(v1);
   \draw[thick] (x1)--(x2)--(x3)--(x4)--(x5)--(x1);
   \draw[thick] (x2)--(x5)--(x3);
   \draw[thick] (y1)--(y2)--(y3)--(y4)--(y1)--(y5)--(y4);
   \draw[thick] (z3)--(z2)--(z1)--(z3)--(z4);

   \end{tikzpicture}
\caption{}\label{fig:special graphs}
\end{figure}


In this paper, we characterized the partitionable (bull, diamond)-free graphs. Moreover, for such graphs, a partition can be found in efficient time.

The next application of the structural behavior of (bull, diamond)-free graphs is on the cops and robber game. The \emph{cops and robber} game on a connected graph is a two-team game. The first team is the team of cops and the second team has a single player named Robber. In the first round, each of the cops chooses a vertex. Then the robber chooses a vertex. Then each team plays in the alternative turn starting with the cops. In the cops' turn, each cop has to choose either to stay at the current place or it can move to an adjacent vertex. The robber also gets the same options to choose (either stays at its current place or moves to an adjacent vertex) in its turn. We say a cop \emph{captured} the robber if the cop is at the same vertex as the robber. The \emph{cop wins} the game if at some turn, one of the cops captures the robber. 
The \emph{cop number} of a connected graph $G$, denoted by $cop(G)$ is the minimum number of cops that can ensure the capture of the robber.  

The first version of the cops and robber game was played with exactly one cop and one robber. This was introduced by Quilliot \cite{quilliot1978} and independently by Nowakowski and Winkler \cite{nowakowski1983}. Later, Aigner and Fromme \cite{aigner1984} introduced the version with multiple cops. Determining the cop number of a graph is known to be $NP$-hard. Therefore the cop number of several graph classes have been studied in the past. 

A family of graphs $\mathcal{F}$ is said to be \emph{cop bounded} if there exists an integer $k$ such that the cop number of any graph in the family is at most $k$. Joret, Kaminski and Theis \cite{joret2008} showed that the class of $H$-free graph is cop-bounded if and only if each component of $H$ is a path. Moreover, they showed the cop number of a connected $P_n$-free graph is at most $n-2$, for any natural number $n>3$. Sivaraman \cite{sivaraman2019application} gave another proof of the same result by using Gy\'arf\'as path argument. In the same paper, he conjectured the following.

\begin{conj}\cite{sivaraman2019application}
    The cop number of a connected $P_n$-free graph is at most $n-3$.
\end{conj}

This conjecture was verified for a specific subclass of $P_5$-free graphs, namely $2K_2$-free graphs (the disjoint union of two edges) \cite{turcotte2022cops}. Further Gupta, Mishra and Pradhan verified several subclasses of $P_5$-free graphs \cite{copP5-free}. Recently Chudnovsky, Norin, Seymour and Turcotte \cite{chudnovsky2024} showed the above conjecture is true when $n$ is $5$. In this paper, we show that the cop number of a connected ($P_n$, bull, diamond)-free graph with a triangle is at most $n-3$. Therefore the cop number of a connected (bull, diamond)-free graph is at most two less than its diameter. Liu \cite{liu2019cop} claimed that the cop number of a ($P_n$, triangle)-free graph is at most $n-3$. So if that is true then the cop number of a connected ($P_n$, bull, diamond)-free is at most $n-3$.

Further, the structural theorems for (bull, diamond)-free graphs given in this paper can have several other applications. A couple of them have been discussed in the conclusion.





\section{Preliminaries}

The vertex set and edge set of a graph $G$ are denoted by $V(G)$ and $E(G)$, respectively. The set of \emph{neighbours} of a vertex $u$ denoted by $N(u)$ is the collection of all the vertices adjacent to $u$ and the closed neighbour of $u$ denoted by $N[u]$ is $\{u\} \cup N(u)$. The \emph{neighbours of a vertex set $S$} denoted by $N(S)$ is the set $\{x \in V(G) \setminus S \mid x \text{ is adjacent to a vertex in } S\}$.  The subgraph \emph{induced} by $S$ in $G$ is denoted by $G[S]$.

The length of a path is the number of edges in the path. The distance from a vertex $u$ to another vertex $v$ in a graph $G$, denoted by $d(u,v)$ is the length of the shortest path from $u$ to $v$. The distance from a vertex $u$ to a subgraph $S$ of $G$, denoted by $d(u, S)$ is defined by $\min_{x \in S} d(x,u)$.

Let $G$ be a graph and $S = \{v_1,v_2,\dots v_t\}$ subset of the vertex set of $G$. We define the following sets for any natural numbers $d$ and $1 \leq i < j \leq t$.

$$W_S(v_i):= \{x \in V \setminus S \mid N(x) \cap S = \{v_i\}\}$$
$$W_S(v_iv_j) := \{v \in V(G) \mid N(v) \cap S = \{v_i,v_j\}\}$$
$$N_d(S) := \{x \in V(G) \mid d(x,S)=d\}.$$

When the context is clear, for every natural numbers $d$ and $1 \leq i < j \leq t$, we write  $W_i, W_{ij}$ and $N_d$ instead of $W_S(v_i)$ and $W_S(v_iv_j)$ and $N_d(S)$. For convenience, we define $N_0 := S$. 

Let $S'$ be another subset of the vertex set. We define $G[S, S']$ to be the graph obtained by collecting all the edges that have one endpoint in $S$ and the other endpoint in $S'$ (the vertex set is the collection of the endpoints of all such edges). 
We say $G[S, S']$ is complete if each vertex in $S$ is a neighbour of every vertex in $S'$. 
A \emph{component} of a graph $G$ is a maximal \textit{connected} subgraph of $G$.
Now we mention some observations on diamond-free graphs.

\begin{observation}\label{prop:NeighTriangle}
Let  $G$ be a diamond-free graph. If $v$ is adjacent to at least two vertices in a clique $K$ then it is adjacent to all the vertices in $K$.
\end{observation}

\begin{observation}\label{prop:Nv}
Let $G$ be a diamond-free graph. For every vertex $v$, $N(v)$ induces a $P_3$-free graph. In other words, $G[N(v)]$ is a disjoint union of cliques.
\end{observation}

\begin{definition}
Let $G$ be a graph. We call $G$ a \emph{matched complete graph} if the vertex set of $G$ can be partitioned into $(A, B)$ such that the following statements hold.

\begin{enumerate}
    \item $G[A]$ and $G[B]$ are complete graphs.
    \item $G[A,B]$ is a disjoint union of edges.
\end{enumerate}
\end{definition}

\section{Structure of (bull, diamond)-free graphs}

We know that complete graphs and matched complete graphs are (bull, diamond)-free. In this section, we examine the structural behavior of other (bull, diamond)-free graphs.

$\mathcal{G} :=$ the collection of all connected (bull, diamond)-free graphs that are neither a complete graph nor a matched complete graph.

Note that the following lemma immediately follows from Observation~\ref{prop:NeighTriangle}.
\begin{lemma}\label{lem:N1UnionOfWis}
    Let $G$ be a connected diamond-free graph and $K = \{v_1,v_2, \dots v_t\}$ be a maximal clique of order at least three. If $G$ is not complete then $(W_1, W_2, \dots W_t)$ is a partition of $N_1$.
\end{lemma}



Note that if exactly one of the $W_i$ is non-empty, then by the above Lemma and Observation~\ref{prop:Nv}, we can conclude that $N_1$ induces a disjoint union of cliques. The following lemma says that $N_1$ induces a complete multipartite graph if at least two of the vertices in a maximal clique have neighbours outside the clique.  

\begin{lemma}\label{lem:MultipartitionOfN1}
Let $G$ be a graph in $\mathcal{G}$ and $K = \{v_1,v_2, \dots v_t\}$ be a maximal clique of order at least three. Suppose $W_i$ and $W_j$ are non-empty sets, for some $1 \leq i< j \leq t$. Then $G[N_1]$ is a complete multipartite graph with partition $(W_1, W_2, \dots W_t)$.
\end{lemma}

\begin{proof}

By Lemma~\ref{lem:N1UnionOfWis}, we know that $(W_1, W_2, \dots W_t)$ is a partition of $N_1$. Without loss of generality, we may assume that $i$ and $j$ are $1$ and $2$, respectively. Let $x$ and $y$ be two vertices in $W_1$ and $W_2$, respectively. Now we claim the following.

\begin{claim}\label{prop:[Wi,Wj]Complete}
    $G[W_k,W_\ell]$ is complete, for any $1 \leq k < \ell \leq t$.
\end{claim}

The set $\{v_1,v_2,v_3,x,y\}$ induces a bull unless $x$ is adjacent to $y$. Since $G$ is bull-free, $xy$ is an edge in $G$. Thus, $G[W_1,W_2]$ is complete. By a similar argument, one can check $G[W_k, W_\ell]$ is complete, for any $1 \leq k < \ell \leq t$. $\Diamond$

Hence, to show $G[N_1]$ is a complete multipartite graph with partition $(W_1, W_2, \dots W_t)$, it is enough to show $W_\ell$ is an independent set, for every $1 \leq \ell \leq t$. Let $u$ and $v$ be two vertices in $W_\ell$, for some $2 \leq \ell \leq t$. By the Claim~\ref{prop:[Wi,Wj]Complete}, we know that $x$ is adjacent to both $u$ and $v$. Now $\{x,u,v,v_\ell\}$ induces a diamond if $u$ is adjacent to $v$. Since $G$ is diamond free, $u$ is not a neighbour of $v$. Hence, $W_\ell$ is an independent set, for any $2 \leq \ell \leq t$. By a similar argument, one can show that $W_1$ is also an independent set. 
\end{proof}

Now we study another important behavior of the structure of neighbours of a maximal clique in a graph in $\mathcal{G}$. In the next lemma, we show that at most two vertices in a maximal clique can have neighbours outside the clique. 

\begin{lemma}\label{lem:W>3empty}
    Let $G$ be a graph in $\mathcal{G}$ and $K = \{v_1,v_2, \dots v_t\}$ be a maximal clique of size at least three. Then there exist $1 \leq i< j \leq t$ such that $W_\ell = \emptyset$, for all $\ell \in \{1,2,\dots, t \} \setminus \{i,j\}$.
\end{lemma}

\begin{proof}

Suppose for the contradiction $W_i, W_j$ and $W_\ell$ are non-empty sets, for some $1 \leq i < j < k \leq t$. Without loss of generality, we may assume that $i,j$ and $k$ are $1,2$ and $3$, respectively. 
By Lemma~\ref{lem:MultipartitionOfN1}, we know that $G[N_i]$ is a complete multipartite graph with partition $(W_1, W_2, \dots W_t)$. In addition, $G$ is a diamond-free graph. As $W_2$ and $W_3$ are non-empty, $W_1$ has exactly one vertex. Similarly, one can argue that, $\mid W_i \mid \leq 1$, for all $1 \leq i \leq t$. Moreover, $N_1$ is a clique. We have assumed that $G$ is a connected graph which is not a matched complete graph. Thus, $N_2$ is a non-empty set.
Now we claim each vertex in $N_1 \cup N_2$ is a clique. 

\begin{claim}\label{G[N1UnionN2]IsComplete}
    $N_1 \cup N_2$ is a clique.
\end{claim}

Let $u$ be a vertex in $N_2$. The definition of $N_2$ says that $u$ has a neighbour in $N_1$. Without loss of generality, we may assume that $u$ is adjacent to the vertex $x$ in $W_1$. Let $y$ and $z$ be the vertices in $W_2$ and $W_3$, respectively. We have just shown that $N_1$ is a clique.
Then $\{u,x,y,z,v_3\}$ induces a bull if $u$ has no neighbour in $\{y,z\}$. Thus, either $y$ or $z$ is a neighbour of $u$. Thus by Observation~\ref{prop:NeighTriangle}, we know that each vertex in $N_1$ is adjacent to $u$. In other words, $G[N_1,N_2]$ is complete. 
Let $v$ be another vertex in $N_2$. Then $\{u,v,x,y\}$ induces a diamond, unless $u$ is adjacent to $v$. Thus, $u$ is a neighbour of $v$. Hence, the above claim holds. $\Diamond$.

We have assumed that $G$ is a connected graph which is not a matched complete graph. Thus, $N_3$ is a non-emptyset. Let $u'$ be a vertex in $N_3$. The definition of $N_3$ says that $u'$ has a neighbour in $N_2$, say $u$. Then $\{v_1,x,y,u,u'\}$ induces a bull. This is a contradiction. 
\end{proof}

Now we are ready to get the complete structure of neighbours of a maximal clique in a graph in $\mathcal{G}$, by using the above two lemmas. 

\begin{theorem}\label{thm:strN1}
Let $G$ be a graph in $\mathcal{G}$ and $K = \{v_1,v_2, \dots v_t\}$ induces a maximal clique of order at least three. Then exactly one of the following statements holds.

\begin{enumerate}
    \item There exists $1 \leq i \leq t$ such that $N_1 = W_i$ and $G[N_1]$ is a disjoint union of cliques.
    
    \item There exist $1 \leq i< j \leq t$ such that $G[N_1(K)]$ is a complete bipartite graph with partition $(W_i,W_j)$.
\end{enumerate}
\end{theorem}

\begin{proof}
By Lemma~\ref{lem:N1UnionOfWis}, we know that $(W_1, W_2, \dots W_t)$ is a partition of $N_1$. Suppose $N_1 = W_i$, for some $1 \leq i \leq t$. Then by Observation~\ref{prop:Nv}, $G[N_1]$ is a disjoint union of cliques. So we may assume that $W_i$ and $W_j$ are two non-empty sets for some $1 \leq i < j \leq t$. Without loss of generality, we may assume that $i$ and $j$ are $1$ and $2$, respectively. By Lemma~\ref{lem:W>3empty}, we know that $W_\ell$ is empty, for any $3 \leq \ell \leq t$. Thus, by Lemma~\ref{lem:MultipartitionOfN1}, $G[N_1]$ is a complete multipartite graph with multipartation $(W_1,W_2)$.

\end{proof} 

In the following two theorems, we study the complete structure of graphs in $\mathcal{G}$.

\begin{theorem}\label{thm:G[N>1]}
Let $G$ be a graph in $\mathcal{G}$ and $K = \{v_1,v_2, \dots v_t\}$  be a maximal clique of order at least three. Then for all $d>1$, $G[N_d]$ is the disjoint union of complete graphs and triangle-free graphs.
\end{theorem}

\begin{proof}

A classical result says that a connected paw-free graph is either triangle-free or a complete-multipartite graph \cite{paw-free}. Further $G$ is a diamond-free graph. Thus it is sufficient to show that $N_d$ induces a paw-free graph for all natural numbers $d > 1$. Suppose for the contradiction, for some $d>1$, $G[N_d]$ contains an induced paw, say with vertex set $\{u_1,u_2,u_3,u_4\}$ and edge set $\{u_1u_2,u_2u_3,u_3u_1,u_1u_4\}$. Further, we claim the following.

\begin{claim}\label{clm:NoComNeighOfu1u2u3}
    $u_1,u_2$ and $u_3$ does not have a common neighbour in $N_{d-1}$.
\end{claim} 

For the contradiction, suppose $v \in N_{d-1}$ is a common neighbour of $u_1,u_2$ and $u_3$. Then $v$ is not a neighbour of $u_4$, otherwise $\{v,u_1,u_2,u_4\}$ induces a diamond. The definition of $N_{d-1}$ says that $v$ has a neighbour in $N_{d-2}$, say $v'$. Now $\{v',v,u_2,u_1,u_4\}$ induces a bull. This is a contradiction. Therefore the above claim is true. $\Diamond$

The definition of $N_d$ says that $u_2$ has a neighour in $N_{d-1}$, say $y$. We know that $\{y,u_1,u_2,u_3,u_4\}$ does not induce a bull. Thus, $y$ has at least one neighbour in $\{u_1,u_3,u_4\}$. By Observation~\ref{prop:NeighTriangle}, $yu_1$ is an edge if and only if $yu_3$ is an edge. Now by claim~\ref{clm:NoComNeighOfu1u2u3}, we get neither $u_1$ nor $u_3$ is a neighbour of $y$. Hence $y$ must be adjacent to $u_4$. In other words, $N(y) \cap \{u_1,u_2,u_3,u_4\} = \{u_2,u_4\}$. By a similar argument, we can show the existence of a vertex $z\in N_{i-1}$ such that $N(z)$ $\cap \{u_1,u_2,u_3,u_4\} = \{u_3,u_4\}$.
Further, $y$ must be adjacent to $z$, as otherwise $\{u_1,u_2,u_3,y,z\}$ induces a bull. 

Again $u_1$ must has a neighbour in $N_{d-1}$, say $x$. By Observation~\ref{prop:NeighTriangle} and claim~\ref{clm:NoComNeighOfu1u2u3}, one can ensures that neither $u_2$ nor $u_3$ is a neighbour of $x$. If $x$ is adjacent to $u_4$ then $\{x,v,u_1,u_2,u_4\}$ induces a bull, for some $v \in N_{d-2} \cap N(x)$. Therefore, $x$ is not a neighbour of $u_4$. However, $y$ must be a neighbour of $x$, otherwise $\{x,y,u_1,u_2,u_3\}$ induces a bull. Similarly, we can show that $xz$ is an edge in $G$.
Hence the graph induced by $\{x,y,z,u_1,u_2,u_3,u_4\}$ is the graph mentioned in Figure \ref{fig:xyzu1u2u3u4}.

\begin{figure}[h!]
  \centering
  \begin{tikzpicture}[scale=1.5]
   \tikzstyle{vertex}=[circle,draw,thick,inner sep=0.2mm,minimum size=1.5mm]

   \node[] (u1) at (0,0) [vertex] {$u_1$};
   \node[] (u2) at (-1,1) [vertex] {$u_2$};
   \node[] (u3) at (1,1) [vertex] {$u_3$};
   \node[] (u4) at (0.5,-1) [vertex] {$u_4$};

   \node[] (x) at (-0.5,-1) [vertex] {$x$};
   \node[] (y) at (-1,-2) [vertex] {$y$};
   \node[] (z) at (1,-2) [vertex] {$z$};

   \draw[thick] (u1) --(u2)--(u3)--(u1)--(u4);
   \draw[thick] (x) -- (y)--(z)--(x);
   \draw[thick] (x)--(u1);
   \draw[thick] (u4)--(y)--(u2);
   \draw[thick] (u4)--(z)--(u3);
   \end{tikzpicture}
\caption{}\label{fig:xyzu1u2u3u4}
\end{figure}

Then the graph induced by $\{x,y,z,u_4\}$ is a diamond. This is a contradiction. Therefore, $N_d$ induces a paw-free graph for all natural number $d>1$.

\end{proof}

\begin{theorem}\label{thm:TriangleinNdNd+1ToNd+2}
Let $G$ be a graph in $\mathcal{G}$ and $K = \{v_1,v_2, \dots v_t\}$ be a maximal clique of order at least three. Then for each $d \geq 1, N(\{x,y,z\}) \cap N_{d+2}$ is an emptyset, whenever $x-y-z-x$ is a triangle in $G[N_d \cup N_{d+1}]$.
\end{theorem}

\begin{proof}
Let $x-y-z-x$ be a triangle in $G[N_d \cup N_{d+1}]$ for some $d \geq 1$. Suppose for the contradiction that $N(\{x,y,z\}) \cap N_{d+2})$ is not an emptyset. Without loss of generality, we may assume that $x \in N_{d+1}$ and it has a neighbour $\Tilde{x}$ in $N_{d+2}$. Now we break our proof into three exclusive cases.

\begin{itemize}
    \item[Case 1] Either $y$ or $z$ is in $N_{d}$
\end{itemize}

Without loss of generality, we may assume that $y$ is in $N_{d}$. Let $y'$ be a neighbour of $y$ in $N_{d-1}$. Then by Observation~\ref{prop:NeighTriangle} one can conclude that either $\{y',x,y,z\}$ induces a diamond or $\{y',x,y,z,x'\}$ induces a bull.
This is a contradiction.



\begin{itemize}
    \item[Case 2] $y,z \in N_{d+1}$ and $\Tilde{x}$ is adjacent to $y$.
\end{itemize}

By Observation~\ref{prop:NeighTriangle}, we know that $\Tilde{x}$ is a neighbour of $z$.
Let $y' \in N_d$ be a neighbour of $y$. We know that $y'$ is not a neighbour of $\Tilde{x}$. By Observation~\ref{prop:NeighTriangle}, neither $x$ nor $z$ is a neighbour of $y'$. So $N(y') \cap \{x,y,z\} = \{y\}$. Similarly, there exists $x',z' \in N_d$ such that $N(x') \cap \{x,y,z\} = \{x\}$ and $N(z') \cap \{x,y,z\} = \{z\}$. Again, $\{x',x,y,z,z'\}$ induces a bull, if $x'$ is not adjacent to $z$. Therefore $x$ must be a neighbour of $z'$. Similarly, we can show that $x'y'$ and $y'z'$ are two edges. The graph induced by $\{\Tilde{x},x,y,z,x',y',z'\}$ is given in Figure\ref{fig:case3}.

\begin{figure}[h!]
  \centering
  \begin{tikzpicture}[scale=0.4]
   \tikzstyle{vertex}=[circle,draw,thick,inner sep=0.2mm,minimum size=1.5mm]

   \node[] (x') at (5,0) [vertex] {$\Tilde{x}$};
   
   \node[] (x) at (1,-5) [vertex] {$x$};
   \node[] (y) at (0,0) [vertex] {$y$};
   \node[] (z) at (1,5) [vertex] {$z$};

   \node[] (x1) at (-5,-5) [vertex] {$x'$};
   \node[] (y1) at (-4,0) [vertex] {$y'$};
   \node[] (z1) at (-5,5) [vertex] {$z'$};

   \draw[thick] (x) --(y)--(z)--(x)--(x');
   \draw[thick] (y) -- (x')--(z);
   \draw[thick] (x1)--(y1)--(z1)--(x1);
   \draw[thick] (x)--(x1);
   \draw[thick] (y)--(y1);
   \draw[thick] (z)--(z1);
   \end{tikzpicture}
\caption{}\label{fig:case3}
\end{figure}

Note that $x'$ has a neighbour in $N_{d-1}$, say $x''$.
Note that $\{x'',x,y,z,x\}$ induces a bull unless $x''$ has a neighbour in $\{z',y'\}$. Without loss of generality, we may assume that $x''$ is adjacent to $y'$.
Observation~\ref{prop:NeighTriangle} says that $x''$ is also adjacent to $z'$. In other words, $x'-y'-z'-x'$ is a triangle. Suppose $d$ is $1$. By Theorem~\ref{thm:strN1}, we know that $N_1 = W_j$, for some $1 \leq j \leq \omega$. Without loss of generality, we may assume that $j=1$ ($x''$ is $ v_1$). Then $\{v_2,x'',x',y',y\}$ induces a bull. This gives a contradiction. Thus, $d$ is at least two. Then for any $v \in ~N(x'') \cap N_{d-2}$, $\{v,x'',x',y',y\}$ induces a bull. This is a contradiction.

\begin{itemize}
    \item[Case 3] $y,z \in N_{d+1}$ and $\Tilde{x}$ is not adjacent to $y$.
\end{itemize}

By Observation~\ref{prop:NeighTriangle}, we know that $\Tilde{x}$ is not adjacent to $z$.
The definition of $N_{d+1}$ says that $y$ must have a neighbour in $N_d$ say $y'$. If $N(y') \cap \{x,z\} = \emptyset$ then $\{y',x,y,z,\Tilde{x}\}$ induces a bull. Thus either $x$ or $z$ is a neighbour of $y'$. By Fact~\ref{prop:NeighTriangle}, $y'$ is adjacent to both $x$ and $z$. We know that $y'$ has a neighbour in $N_{d-1}$, say $y''$. Then $\{y'',y',y,x,\Tilde{x}\}$ induces a bull. This gives a contradiction.

\end{proof}


\section{Partitioning (bull, diamond)-free graphs}

In this section, we characterize the (bull, diamond)-free partitionable graphs. Further, we show that for partitionable (bull, diamond)-free graphs, we can find such a partition in polynomial time.




\begin{theorem}\label{thm:NotInMathcalGPartition}
    Let $G$ be a graph in $\mathcal{G}$ that has a triangle. Then $G$ is partitionable. Moreover, such a partition can be found in polynomial time.
\end{theorem}

\begin{proof}
First, we find a maximal clique containing a size at least three, say $K = \{v_1,v_2,v_3,\dots t\}$ (this step takes $O(n^4)$-time). We can find the partition of $(K,N_1,N_2, \dots N_\ell)$ of $G$ in polynomial time, where $\ell$ is the largest distance of a vertex from $K$. Since $G$ is connected but not a complete graph, without loss of generality, we may assume that $W_1$ has a vertex, say $w_1$. By Theorem~\ref{thm:G[N>1]} we know that $1 \leq d \leq t$ for any $1 \leq d \leq \ell$, $G[N_d]$ is disjoint union of cliques and triangle-free graphs. Now we define the following sets (that can be found in polynomial time) for each $1 \leq d \leq \ell$, which is helpful in getting a partition of $G$.

\begin{align*}
    A_d := \{x \in N_d \mid & \, x \text{ is a vertex in a clique of size at least three in } G[N_d]\}\\
    B_d :=  \{x,y \in N_d \mid & \, \exists z \in N_{d-1} \text{ such that } x-y-z-x \text{ is a triangle and } \\ & \{x,y\} \text{ has no neighbour in } N_d\}\\
    C_d :=  \{x \in N_d \mid & \, x \text{ is adjacent to no vertex in $N_d$ and } \\ & x-y-z-x \text{ is a triangle, for some } y,z \in N_{d-1}\}\\
\end{align*}

Note that for any $d>1$, $A_d,B_d$ and $C_d$ mutually non-intersecting sets. Further, any vertex in $C_d$ is an isolated vertex in $G[N_d]$. Further $G[B_d]$ is a disjoint union edges and $G[A_d]$ is a disjoint union of cliques.  Moreover, no vertex in $B_d$ has a neighbour in $A_d$. Thus, $A_d \cup B_d \cup C_d$ induces a disjoint union of cliques, for any $1 < d \leq t$. By Theorem~\ref{thm:TriangleinNdNd+1ToNd+2} we know that no vertex in $A_d \cup B_d \cup C_d$ has a neighbour in $N_{d+1}$, for all $1 < d \leq t$. Therefore, the following set induces a disjoint union of cliques.

$$\mathcal{A} := \bigcup_{1 < d \leq t} A_d \cup \bigcup_{1 < d \leq t} B_d \cup \bigcup_{1 < d \leq t} C_d$$

Thus, $G \setminus (K \cup N_1)$ is partitionable with partition $(\mathcal{A}, G \setminus (K \cup N_1 \cup \mathcal{A}))$. Now we divide our proof into two cases depending upon whether $W_i = \emptyset$, for all $1 < i \leq t$. 

\begin{itemize}
    \item[Case 1] $W_i = \emptyset$, for all $1 < i \leq t$
\end{itemize}

Thus, $N_1 =W_1$ and hence $G[N_1]$ induces disjoint union of cliques. Lemma~\ref{thm:TriangleinNdNd+1ToNd+2} ensures us that $G$ is partitionable with partition $(\mathcal{A} \cup \mathcal{B} \cup \{v_2,v_3, \dots v_t\}, V(G) \setminus (\mathcal{A} \cup \mathcal{B} \cup \{v_2,v_3, \dots v_t\}))$ is a partition of $G$, where $\mathcal{B}$ is the set of all non-isolated vertices in $G[W_1]$.

\begin{itemize}
    \item[Case 2] $W_i \neq \emptyset$, for some $1 < i \leq t$
\end{itemize}

Without loss of generality, we may assume that $W_2$ is a non-empty set. By Lemma~\ref{lem:W>3empty}, we know that $W_i = \emptyset$, for all $3 \leq i \leq t$ and $G[N_1]$ is a complete bipartite graph with bipartition $(W_1,W_2)$. 
Then $G$ is partitionable with partition $(\mathcal{A} \cup \{v_2,v_3,v_4, \dots v_t\}, V(G) \setminus (\mathcal{A} \cup \{v_2,v_3,v_4, \dots v_t\}))$.





\end{proof}

\begin{theorem}
There is a polynomial-time algorithm that decides whether a connected (bull, diamond)-free graph is partitionable. Further when it is partitionable such a partition can be found in polynomial-time.
\end{theorem}

\begin{proof}
Let $G$ be the input graph. Note that if $G$ is triangle-free, then it is partitionable with the obvious partition.
Hence, first, we check whether $G$ has a triangle. Moreover if yes then we get a triangle $v_1-v_2-v_3-v_1$ (the running time for this is $O(n^3)$-time). Moreover, we find a maximal clique $K$ containing $v_1,v_2$ and $v_3$. Again, if $G$ is a complete graph then it is partitionable with partition $(V(G), \emptyset)$. If $G$ is a matched complete graph then it is partitionable if and only if $G[K, V(G) \setminus K]$ is of size at most two. Further, in this case, the $G$ is partitionable with partition $(K, V(G) \setminus K)$. So we may assume that $G \in \mathcal{G}$. Now we follow Theorem~\ref{thm:NotInMathcalGPartition}.




\end{proof}





    

    

\section{Cop number of (bull, diamond)-free graphs}

In this section, we discuss the cops and robber game on (bull, diamond)-free graphs with triangles. The proof of the following theorem is inspired by the Gy\'arf\'as path argument \cite{gyarfas}.

\begin{theorem}\label{thm:CopNoGraphsWithTriangle}
Let $G$ be a connected ($P_n$, bull, diamond)-free graph. If $G$ has a triangle then the $c(G) \leq n-3$.
\end{theorem}

\begin{proof}

Let $K = \{v_1, v_2, \dots v_t\}$ be a maximal clique of size at least three in $G$. The cop number of a complete graph is $1$ and that of a matched complete graph is $2$. So we may assume that $G$ is in $\mathcal{G}$. So $N_1$ is a nonempty set. Lemma~\ref{lem:N1UnionOfWis} says that $(W_1, W_2, \dots W_t)$ is a partition of the vertex set of $N_1$. Without loss of generality, we may assume that $W_1$ is a non-empty set. Now we split our proof into two cases.

\begin{itemize}
    \item[Case 1] $W_i = \emptyset$, for all $1 < i \leq t$
\end{itemize}

Now we play the game. In the first round, we place all the $n-3$ cops at $v_1$. For convenience, we write $u_0=v_1$. Now the robber has to choose a vertex in $G \setminus N[u_0] = N_2 \cup N_3 \cup \dots \cup N_{n-1}$. Let $C_1$ be the component of $G \setminus N[u_0]$ that contains the robber and $u_1$ be a neighbour of $u_0$ that has a neighbour in $C_1$. In the next round one cop stays in $u_0$ and the rest $n-4$ robbers move to $u_1$. In the next round robber has to be in a component of $C_1 \setminus N[u_1]$, say $C_2$. Let $u_2$ be a neighbour of $u_1$ that has a neighbour in $C_2$. In the next step; the cop at $u_0$ stays there, one cop stays at $u_1$ and the rest of the cops move to $u_2$. We repeat this process $n-5$ times. Now we have a nested sequence of induced subgraphs $C_1 \supset C_2 \supset C_3 \dots C_{n-4}$ and an induced path $u_0-u_1-u_2-\dots-u_{n-4}$ in which every vertex has a cop. 

Now we claim that in the next step, the robber can move only to a neighbour of $u_i$, for some $0 \leq i \leq n-4$. Suppose for the contradiction the robber moves to a vertex $r$ which does not has a neighbour in $\{u_0, u_1 \dots u_{n-4}\}$. Let $P$ be the shortest path from a neighbour of $u_{n-4}$ to $r$ in $C_{n-4}$. Note that $N_1=W_1$. So $v_2$ has no neighbour in $\{u_1,u_2, \dots, u_{n-4}\}$. Thus the vertices of $P$ together with $v_2$ and $\{u_0, u_1 \dots u_{n-4}\}$ gives a path of length at least $n$. This is a contradiction. Therefore the robber has to choose a neighbour of $u_i$, for some $0 \leq i \leq n-4$. Hence, it gets captured in the next round. 

\begin{itemize}
    \item[Case 2] $W_i \neq \emptyset$, for some $1 < i \leq t$
\end{itemize}

Without loss of generality, we may assume that $W_2$ is a non-empty set. By Lemma~\ref{lem:W>3empty} we know that $W_i = \emptyset$, for all $3 \leq i \leq t$. Further by Lemma~\ref{thm:strN1} we know that $G[N_1]$ is a complete bipartite graph with bipartition $(W_1,W_2)$.

Now we give a winning strategy for the cops. In the first round, we place all the $n-3$ cops in $v_1$. Then robber has to choose a vertex in $G \setminus (K \cup W_1)$. Let $C_1$ be the component of $G \setminus (K \cup W_1)$ that contains the robber. Note that one of the neighbour of $v_1$ has a neighbour in $C_1$. Since $N_1 = W_1 \cup W_2$, either $v_2$ or a vertex in $W_1$ has a neighbour in $C_1$. Moreover, if $v_2$ has a neighbor in $C_1$. Then one of the vertex of $W_2$ is in $C_1$. Theorem~\ref{thm:strN1} says that $G[N_1]$ is a complete bipartite graph with bipartition $(W_1,W_2)$. So each vertex in $W_1$ has a neighbour in $C_1$. Therefore, we may assume that one of the vertex of $W_1$ has a neighbour in $C_1$, say $u_1$.

In the next round, one of the cops stays at $v_1$ and the rest of the cops move to $u_1$. Now we follow a strategy similar to the previous case. In the next round robber has to be in a component of $C_1 \setminus N[u_1]$, say $C_2$ (Note that $C_2$ has no vertex in $N_1$). Let $u_2$ be a neighbour of $u_1$ that has a neighbour in $C_2$.  In the next step; the cop at $v_1$ stays there, one cop stays at $u_1$ and the rest of the cops move to $u_2$. We repeat this process $n-5$ times. So now we may assume that we have an induced path $v_1-u_1-u_2-\dots-u_{n-4}$ and each of the vertex in this path has one cop. Additionally, we have a nested sequence of induced subgraphs $C_1 \supset C_2 \supset C_3 \dots C_{n-4}$.

Now we claim that in the next step, the robber can move only to a neighbour of $\{v_1, u_1 \dots u_{n-4}\}$. Suppose for the contradiction the robber moves to a vertex $r$ which does not have a neighbour in $\{v_1, u_1 \dots u_{n-4}\}$. Let $P$ be the shortest path from $u_3$ to $r$ in $C_{n-4}$. Then the vertices of $P$ together with $v_3$ and $\{v_1, u_1 \dots u_{n-4}\}$ gives a path of length at least $n$. This is a contradiction. Therefore the robber has to choose a neighbour of $\{v_1, u_1 \dots u_{n-4}\}$. Hence, it gets captured in the next round. 

\end{proof}

\begin{corollary}
The cop number of a (bull, diamond)-free graph with a triangle is at most two less than its diameter.
\end{corollary}

Liu \cite{liu2019cop} claimed that the cop number of a connected ($P_n$, triangle)-free graph is at most $n-3$. 
If that is true, then by Theorem~\ref{thm:CopNoGraphsWithTriangle} the cop number of a connected ($P_n$, bull, diamond)-free graph is at most $n-3$.

\section{Conclusion}

The structural behaviors of the (bull, diamond)-free graphs which we study in Section $3$ have multiple applications. Here, we discuss two applications of Lemma~\ref{lem:W>3empty}: clique covering number and clique chromatic number of (bull, diamond)-free graphs. 

The clique covering number of a graph $G$ denoted by $\theta (G)$ is the minimum number of cliques in $G$ needed to cover the vertex set of $G$. The clique covering number of a matched complete graph is two and that of a complete graph is one. To under the clique covering number of (bull, diamond)-free graphs in $\mathcal{G}$, let us represent Lemma\ref{lem:W>3empty} in another way and for that first we give the following definition. 


\begin{definition}
    Let $G$ be a graph and $x$ and $y$ be two vertices. We say $y$ is \emph{dominated} by $y$ if $N(y) \cup \{y\} \subset N(x) \cup \{x\}$.
\end{definition}

Lemma~\ref{lem:W>3empty} says that every (bull, diamond)-free graph $G$ with a triangle has a dominating vertex. Further, $\theta (G) = \theta (G \setminus \{x\})$, where $G$ is a graph and $x$ is a dominating vertex. Thus, the following corollary follows. 



\begin{corollary}
Let $G \in \mathcal{G}$ be a connected graph. Suppose $\Tilde{G}$ be the graph obtained from $G$ by repeatedly deleting dominated vertices. If $\Tilde{G}$ does not have a dominated vertices then $\theta (G) = \mid V(\Tilde{G}) \mid - m(\Tilde{G})$, where $m(\Tilde{G})$ is the size of a maximum matching of $\Tilde{G}$.
\end{corollary}


One of the other applications is clique-coloring. A \emph{$c$-clique-coloring} of a graph $G$ is an assignment of $c$ colors to the vertices of $G$ such that no maximal clique of size at least $2$ in $G$ is monochromatic. The \emph{clique chromatic number} of $G$ is the minimum number of colors used in a clique coloring of $G$.

Let $G$ be a connected (bull, diamond)-free graph. We have just discussed that if $G$ is not a triangle-free graph. Then $G$ has a dominating vertex. Therefore the clique chromatic number of $G$ is $\chi(G')$, where $G'$ is the graph obtained from $G$ by deleting the dominating vertices in $G$. We know that the chromatic number of a ($P_n$-triangle)-free graph is at most $n-2$ \cite{gravier}. Therefore the clique chromatic number of $G$ is at most $n-2$ where $n$ is the diameter of $G$. 


\section*{Acknowledgement:} This research work is partially funded by Fondecyt Postdoctoral grant $3220618$ of Agencia National de Investigati\'{o}n y Desarrollo (ANID), Chile.

\bibliographystyle{plain}

\end{document}